\documentclass{article}

\usepackage{amssymb}
\usepackage{amsmath}

\newtheorem{theorem}{Theorem}
\newtheorem{lemma}{Lemma}[section]
\newtheorem{propo}{Proposition}[section]

\newtheorem{defin}{Definition}[section]

\newcommand{\N}{\mathbb{N}}

\newcommand{\cT}{\mathcal{T}}
\newcommand{\cF}{\mathcal{F}}

\newcommand{\Grd}[1][d]{\textrm{\textbf{Gr}}_{#1}}
\newcommand{\Gd}[1][d]{\textrm{\textbf{Graph}}_{#1}}

\def\ogr{\vec{\Grd}}
\def\ogrn{\vec{\Grd^n}}
\def\our{\vec{U}^{r,d}}
\def\oa{\vec{\alpha}}\def\ob{\vec{\beta}}
\def\omu{\vec{\mu}}
\def\urd{U^{r,d}}
\def\erd{\vec{E}^{r,d}}
\def\vep{\vec{\phi}}
\def\veps{\vec{\psi}}
\def\we{w_{\delta}}

\def\grnd{\Grd^n}
\def\urnd{\urd_n}
\def\mun{\mu_n}

\newcommand{\defeq}{\stackrel{\mbox{\scriptsize def}}{=}}

\newcommand{\G}{\textrm{\textbf{G}}}

\newcommand{\mug}{\mu_{\bf G}}
\def\limn{\lim_{n\to\infty}}

\def\proof{\smallskip\noindent{\it Proof.} }
\newcommand{\qed} {\hspace {0.1in} \rule {1.5mm} {3.5mm}}
\title{On the limit of large girth graph sequences}
\author{G\'abor Elek}

\begin{document}
\maketitle

\begin{abstract} Let $d\geq 2$ be given and let $\mu$ be an 
involution-invariant probability measure on 
the space of trees $T\in\cT_d$ with maximum
degrees at most $d$. Then $\mu$ arises as the local limit of some sequence
$\{G_n\}^\infty_{n=1}$ of graphs with all degrees at most $d$.
This answers Question 6.8 of Bollob\'as and Riordan \cite{Bol}.
\end{abstract}

\section{Introduction}
Let $\Gd$ denote the set of all finite simple
 graphs $G$ (up to isomorphism)
  for which $\deg(x) \leq d$ for every $x \in V(G)$.
For a graph $G$ and $x,y \in V(G)$ let $d_G(x,y)$ denote the distance
of $x$ and $y$, that is the length of the shortest path from $x$ to
$y$.  A rooted $(r,d)$-ball is a graph $G \in \Gd$ with a marked
vertex $x \in V(G)$ called the root such that $d_G(x,y) \leq r$ for
every $y \in V(G)$.  By $U^{r,d}$ we shall denote the set of rooted
$(r,d)$-balls.

If $G \in \Gd$ is a graph and $x\in V(G)$ then $B_r(x) \in U^{r,d}$
shall denote the rooted $(r,d)$-ball around $x$ in $G$.
For any $\alpha \in U^{r,d}$ and $G \in \Gd$ we define the set
  $T(G,\alpha) \defeq \{x \in V(G): B_r(x) \cong \alpha\}$ and let
  $p_G(\alpha) \defeq \frac{|T(G,\alpha)|}{|V(G)|}$.
A graph sequence $\G = \{G_n\}_{n=1}^{\infty} \subset \Gd$ is
  {\bf weakly convergent} if $\lim_{n \to \infty} |V(G_n)| = \infty$
  and for every $r$ and every $\alpha \in U^{r,d}$ the limit
  $\lim_{n\to\infty}p_{G_n}(\alpha)$ exists (see \cite{BS}).

\noindent
 Let $\Grd$ denote the set of all countable, connected rooted
graphs
  $G$ for which $\deg(x) \leq d$ for every $x \in V(G)$.
If $G,H\in\Grd$ let $d_g(G,H)=2^{-r}$, where
$r$ is the maximal number such that the $r$-balls around the roots
of $G$ resp. $H$ are rooted isomorphic. The distance $d_g$  makes
$\Grd$ a compact metric space. Given an
  $\alpha \in U^{r,d}$ let $T(\Grd,\alpha) = \{(G,x) \in \Grd : B_r(x)
  \cong\alpha\}$. The sets $T(\Grd,\alpha)$ are closed-open sets.
A convergent graphs sequence $\{G_n\}^\infty_{n=1}$ define a 
{\bf local limit measure} $\mu_{\bf G}$ on $\Grd$,
where $\mug(T(\Grd,\alpha)) =\limn p_{G_n}(\alpha)$.
However, not all the probability measures on $\Grd$ arise as local
limits. A necessary condition for a measure $\mu$ being a local limit
is its {\bf involution invariance} (see Section \ref{invol}).
The goal of this paper is to answer a  question of
Bollob\'as and Riordan (Question 6.8 \cite{Bol}):
\begin{theorem}\label{theorem1}
Any involution-invariant measure $\mu$ on $\Grd$ concentrated on
trees arises as a local limit of some convergent graph sequence.
\end{theorem}
As it was pointed out in \cite{Bol} 
such graph sequences are asymptotically treelike, thus
$\mu$ must arise as the local limit of a convergent large girth
sequence.
\section{Involution invariance} \label{invol}
Let $\ogr$ be the compact space of all connected
countable rooted graphs  $\vec{G}$ (up to isomorphism)
of vertex degree bound $d$
with a distinguished directed edge pointing out from the root.
Note that $\vec{G}$ and $\vec{H}$ are considered isomorphic if there
exists a rooted isomorphism between them mapping distinguished
edges into each other.
Let $\our$ be the isomorphism classes
of all rooted $(r,d)$-graphs $\oa$ with a distinguished edge
$e(\oa)$ pointing out from the root.
Again, $T(\ogr,\oa)$ is well-defined for any $\oa\in\our$ and defines
a closed-open set in $\ogr$. Clearly,
the forgetting map $\cF:\ogr\to\Grd$ is continuous.
Let $\mu$ be  a probability measure on $\Grd$. Then we define
a measure $\omu$ on $\ogr$ the following way.

\noindent
Let $\oa\in \our$ and let $\cF(\oa)=\alpha\in\urd$ be the
underlying rooted ball. Clearly,
$\cF(T(\ogr,\oa))=T(\Grd,\alpha)$.Let
$$\omu(T(\ogr,\oa)):=l\,,$$ where
$l$ is the number of edges $e$ pointing out from the root such that there
exists a rooted automorphism of $\alpha$ mapping $e(\oa)$ to $e$.
Observe that
$$\omu(\cF^{-1}(T(\Grd,\alpha))=\deg(\alpha)\mu(T(\Grd,\alpha))\,.$$
We define the map $T:\ogr\to\ogr$ as follows. Let $T(\vec{G})=\vec{H}$,
where :
\begin{itemize}
\item the underlying graphs of $\vec{G}$ and $\vec{H}$ are the same,
\item the root of $\vec{H}$ is the endpoint of $e(\vec{G})$,
\item the distinguished edge of $\vec{H}$ is pointing
to the root of $\vec{G}$.
\end{itemize}
Note that $T$ is a continuous involution.
Following Aldous and Steele \cite{AS}, we
call $\mu$ {\bf involution-invariant} if $T_*(\omu)=\omu$.
It is important to note \cite{AS},\cite{AL} that the limit measure
of convergent graphs sequences are always involution-invariant.

\noindent
We need to introduce the notion of {\bf edge-balls}.
Let $\vec{G}\in\ogr$. The edge-ball $B^e_r(\vec{G})$ of radius
$r$ around the root of $\vec{G}$ is the following spanned rooted
subgraph of $\vec{G}$:
\begin{itemize}
\item The root of $B^e_r(\vec{G})$ is the same as the root of $\vec{G}$.
\item $y$ is a vertex of $B^e_r(\vec{G})$ if
$d(x,y)\leq r$ or $d(x',y)\leq r$, where $x$ is the root of $\vec{G}$
and $x'$ is the endpoint of the directed edge $e(\vec{G})$.
\item The distinguished edge of $B^e_r(\vec{G})$ is $(\vec{x,x'})$.
\end{itemize}
Let $\erd$ be the set of all edge-balls of radius $r$ up to isomorphism.
Then if $\vec{\phi}\in\erd$, let
$s(\vep)\in \our$ be the rooted ball around the root of $\vep$.
Also, let $t(\vep)\in\our$ be the $r$-ball around $x'$ with
distinguished edge $(\vec{x',x})$.

\noindent
The involution $T^{r,d}:\erd\to\erd$ is
defined the obvious way and $t(T^{r,d}(\vep))=s(\vep)$,
$s(T^{r,d}(\vep))=t(\vep)$.
Since $\omu$ is a measure we have 
\begin{equation} \label{e1}
\omu(T(\ogr,\oa))=\sum_{\vep,s(\vep)=\oa}\omu(T(\ogr,\vep)).
\end{equation}
Also, by the involution-invariance
\begin{equation} \label{e3}
\omu(T(\ogr,\vep))=\omu(T(\ogr,T^{r,d}(\vep)),
\end{equation}
since $T(T(\ogr,\vep))=T(\ogr,T^{r,d}(\vep)$.
Therefore by (\ref{e1}),
\begin{equation} \label{e2}
\omu(T(\ogr,\oa))=\sum_{\vep,t(\vep)=\oa}\omu(T(\ogr,\vep))
\end{equation}
\section{Labeled graphs}
Let $\ogrn$ be the isomorphism classes of
\begin{itemize}
\item connected countable rooted graphs with vertex degree bound $d$
\item with a distinguished edge pointing out from the root
\item with vertex labels from the set $\{1,2,\dots,n\}$.
\end{itemize}
Note that if $\vec{G}_*$ and $\vec{H}_*$ are such graphs then
they called isomorphic if there exists a map $\rho:V(\vec{G}_*)\to
V(\vec{H}_*)$ preserving both the underlying $\ogr$-structure and the
the vertex labels.
The labeled $r$-balls $\our_n$ and the labeled $r$-edge-balls
$\erd_n$ are defined accordingly.
Again, $\ogrn$ is a compact metric space and
$T(\ogrn,\oa_*)$,\,$T(\ogrn,\vep_*)$ are closed-open sets,
where $\oa_*\in \our$, $\vep_*\in \erd_n$.
Now let $\mu$ be an involution-invariant probability
measure on $\Grd$ with induced measure $\omu$. The associated
measure $\omu_n$ on $\ogrn$ is defined the following way.

\noindent
Let $\oa\in\our$ and $\kappa_1,\kappa_2$ be vertex labelings of $\oa$ by
$\{1,2,\dots,n\}$. We say that $\kappa_1$ and $\kappa_2$ are equivalent
if there exists a rooted automorphism of $\oa$ preserving the
distinguished edge and mapping $\kappa_1$ to $\kappa_2$.
Let $C(\kappa)$ be the equivalence class of the vertex labeling $\kappa$
of $\oa$. Then we define
$$\omu_n(T(\ogrn,[\kappa])):=\frac{|C(\kappa)|}{n^{|V(\oa)|}}
\omu(T(\ogr,\oa))\,.$$
\begin{lemma} \label{borel}
\begin{enumerate}
\item $\omu_n$ extends to a Borel-measure.
\item $\omu(T(\ogr,\oa))=\sum_{\oa_*,\,\cF(\oa_*)=\oa}
\omu_n(T(\ogrn,\oa_*))\,.$
\end{enumerate} \end{lemma}
\proof
The second equation follows directly from th definition. In order to prove
that $\omu_n$ extends to a Borel-measure it is enough to prove that
$$\omu_n(T(\ogrn,\oa_*))=
\sum_{\ob_*\in N_{r+1}(\oa_*)}
\omu_n(T(\ogrn,\ob_*))\,,$$
where $\oa_*\in\our_n$ and $N_{r+1}(\oa_*)$ is the
set of elements $\ob_*$ in $\vec{U}^{r+1,d}_n$ such that the
$r$-ball around the root of $\ob_*$ is isomorphic to $\oa_*$.
Let $\oa=\cF(\oa_*)\in\our$ and let
$N_{r+1}(\oa)\subset \our$ be the set of elements $\ob$ such that
the $r$-ball around the root of $\ob$ is isomorphic to $\oa$.
Clearly
\begin{equation} \label{gyors}
\omu(T(\ogr,\oa))=\sum_{\ob\in N_{r+1}(\oa)} \omu(T(\ogr,\ob))\,.
\end{equation}
Let $\kappa$ be a labeling of $\oa$ by $\{1,2,\dots,n\}$ representing
$\oa_*$. For $\ob\in N_{r+1}(\oa)$ let $L(\ob)$ be the set of labelings
of $\ob$ that extends some labeling of $\oa$ that is equivalent to $\kappa$.

\noindent
Note that
$$\omu_n(T(\ogrn,\oa_*))=\omu(T(\ogr,\oa))\frac{|C(\kappa)|}{n^{|V(\oa)|}}\,.$$
Also,
$$\sum_{\ob_*\in N_{r+1}(\oa_*)} \omu_n (T(\ogrn,\ob_*))=
\sum_{\ob\in N_{r+1}(\oa)} \omu(T(\ogr,\ob))
\frac{|L(\ob)|}{n^{|V(\ob)|}}\,.$$
Observe that
$|L(\ob)|=|C(\kappa)| n^{|V(\ob)|-V(\oa)|}$.
Hence
$$\sum_{\ob_*\in N_{r+1}(\oa_*)} \omu_n (T(\ogrn,\ob_*))=
\sum_{\ob\in N_{r+1}(\oa)}\omu(T(\ogr,\ob)) 
\frac{|C(\kappa)|}{n^{|V(\oa)|}}\,.$$
Therfore using equation (\ref{gyors}) our lemma follows. \qed

\vskip0.2in
\noindent
The following proposition shall be crucial in our construction.
\begin{propo}\label{master}
For any $\oa_*\in\our_n$ and $\veps_*\in\erd_n$
\begin{itemize}
\item 
$\omu_n(T(\ogrn,\oa_*))=\sum_{\vep_*\in\erd_n,\,s(\vep_*)=\oa_*}
\omu_n(T(\ogrn,\vep_*))$
\item
$\omu_n(T(\ogrn,\oa_*))=\sum_{\vep_*\in\erd_n,\,t(\vep_*)=\oa_*}
\omu_n(T(\ogrn,\vep_*))$
\item
$\omu_n(T(\ogrn,\veps_*))=
\omu_n(T(\ogrn,T^{r,d}_n(\veps_*))\,.$
\end{itemize}
\end{propo}
\proof
The first equation follows from the fact that $\omu_n$ is a Borel-measure.
Thus the second equation will be an immediate corollary of the third one.
So, let us turn to the third equation.
Let $\cF(\veps_*)=\veps\in\erd$ and let $\kappa$ be a vertex-labeling
of $\veps$ representing $\veps_*$. It is enough to prove that
$$\omu_n(T(\ogrn,\veps_*))=
\frac{|C(\kappa)|}{n^{|V(\veps)|}} \omu(T(\ogr,\veps))\,,$$
where $C(\kappa)$ is the set of labelings of $\veps$ equivalent to $\kappa$.
Let $N_{r+1}(\veps)\in\our$ be the set of elements $\ob$ such that the 
edge-ball of radius $r$ around the root of $\ob$ is isomorphic to $\veps$.
Then
\begin{equation} \label{gyors2}
\omu(T(\ogr,\veps))=\sum_{\ob\in N_{r+1}(\veps)} \omu(T(\ogr,\ob))\,.
\end{equation}
Observe that
$$\omu_n(T(\ogrn,\veps_*))=\sum_{\ob\in N_{r+1}(\veps)} \omu(T(\ogr,\ob))
\frac{k(\ob,\veps_*)}{n^{|V(\ob)|}}\,,$$
where $k(\ob,\veps_*)$ is the number of labelings
of $\ob$ extending an element that is equivalent to $\kappa$.
Notice that $k(\ob,\veps_*)=|C(\kappa)|n^{|V(\ob)|-|V(\veps)|}\,.$
Hence by (\ref{gyors2})
$\omu_n(T(\ogrn,\veps_*))=
\frac{|C(\kappa)|}{n^{|V(\veps)|}} \omu(T(\ogr,\veps))\,,$ thus our
proposition follows. \qed
\section{Label-separated balls}
Let $\grnd$ be the isomorphism classes of
\begin{itemize}
\item connected countable rooted graphs with vertex degree bound $d$
\item with vertex labels from the set $\{1,2,\dots,n\}$.
\end{itemize}
Again, we define the space of labeled $r$-balls $\urd_n$. Then
$\grnd$ is a compact space with closed-open sets $T(\grnd,M), M\in\urnd$.
Similarly to the previous section we define an associated probability measure
$\mun$, where $\mu$ in an involution-invariant probability measure on $\Grd$.

\noindent
Let $M\in\urnd$ and let $R(M)$ be the set of  elements of
$\our_n$ with underlying graph $M$.
If $A\in R(M)$, then the multiplicity of $A$, $l_A$ is the number of
edges $e$ pointing out from the root of $A$ such that
there is a label-preserving rooted automorphism of $A$ moving the
distinguished edge to $e$.
Now let
$$\mun(M):=\frac{1}{\deg(M)}\sum_{A\in R(M)} l_A\omu_n(A)\,.$$
The following lemma is the immediate consequence of Lemma \ref{borel}.
\begin{lemma} \label{egyenletek}
$\mu_n$ is a Borel-measure on $\grnd$ and
$\sum_{M\in M(\alpha)}\mu_n(M)=\mu(A)$ if
$\alpha\in\urd$ and $M(\alpha)$ is the set of labelings of $\alpha$ by
$\{1,2,\dots,n\}$.
\end{lemma}
\begin{defin}
$M\in\urnd$ is called label-separated
if all the labels of $M$ are different.
\end{defin}
\begin{lemma}
For any $\alpha\in\urd$ and $\delta>0$ there exists an $n>0$
such that
$$|\sum_{M\in M(\alpha),\,M\,\,\mbox{ is label-separated}}
\mu_n(T(\Grd,M))-\mu(T(\Grd,\alpha))|<\delta\,.$$
\end{lemma}
\proof
Observe that
$$\sum_{M\in M(\alpha),\,M\,\,\mbox{ is label-separated}}
\mu_n(T(\Grd,M))=\frac{T(n,\alpha)}{n^{|V(\alpha)|}}\mu(T(\Grd,\alpha))\,,$$
where $T(n,\alpha)$ is the number of 
$\{1,2,\dots,n\}$-labelings of $\alpha$ with
different labels.
Clearly, $\frac{T(n,\alpha)}{n^{|V(\alpha)|}}\to 1$ as $n\to\infty$. \qed

\section{The proof of Theorem \ref{theorem1}}
Let $\mu$ be an involution-invariant probability measure on $\Grd$
supported on trees.
It is enough to prove that for any $r\geq 1$ and $\epsilon >0$ there exists a
finite graph $G$ such that for any $\alpha\in\urd$
$$|p_G(\alpha)-\mu(T(\Grd,\alpha))|<\epsilon\,.$$
The idea we follow is close to the one used by Bowen in \cite{Bow}.
First, let $n>0$ be a natural number such that
\begin{equation} \label{becs1}
|\sum_{M\in M(\alpha),\,M\,\,\mbox{ is label-separated}}
\mu_n(T(\Grd,M))-\mu(T(\Grd,\alpha))|<\frac{\epsilon}{10}\,.
\end{equation}
\noindent
Then we define a directed labeled finite graph $H$ to encode some
information on $\omu_n$. If $A\in\vec{U}^{r+1,d}_n$ then let $L_A$ be the
unique element of $\erd_n$ contained in $A$.

\noindent
The set of vertices of $H$; $V(H):=\vec{U}^{r+1,d}_n$.
If $A,B\in \vec{U}^{r+1,d}_n$ and $L_A=L^{-1}_B$ (we use
the inverse notation instead of writing out the involution operator) then
there is a directed edge $(A,L_A,B)$ from $A$ to $B$ labeled by $L_A$ and
a directed edge $(B,L_B,A)$ from $B$ to $A$ labeled by $L_B=L^{-1}_A$.
Note that we might have loops.
We define the weight function $w$ on $H$ by
\begin{itemize}
\item $w(A)=\omu_n(T(\ogrn,A))$.
\item $w(A,L_A,B)=\mu(T(\ogrn,L_{A,B}))\,,$
where $L_{A,B}\in\vec{E}^{r+1,d}_n$ the unique element such that
$s(L_{A,B})=A, t(L_{A,B})=B$.
\end{itemize}
By Proposition \ref{master} we have the following equation for
all $A,B$ that are connected in $H$:
\begin{equation} \label{d1}
w(A,L_A,B)=w(B,L^{-1}_A,A)\,.
\end{equation}
Also,
\begin{equation} \label{d2}
w(A)=\sum_{w(A,L_A,B)\in E(H)} w(A,L_A,B)
\end{equation}
\begin{equation} \label{d3}
w(A)=\sum_{w(B,L^{-1}_A,A)\in E(H)} w(B,L^{-1}_A,A)
\end{equation}
Also if $M\in U^{r+1,d}_n$ then
\begin{equation} \label{d4}
\mu_n(M)=\frac{1}{\deg{(M)}}\sum_{A\in R(M)} l_A w(A),
\end{equation}
where $l_A$ is the multiplicity of $w(A)$.

\vskip0.2in
\noindent
Since the equations (\ref{d1}), (\ref{d2}), (\ref{d3}) have rational
coefficients we also have weight functions $\we$ on $H$
\begin{itemize}
\item taking only rational values
\item satisfying equations (\ref{d1}), (\ref{d2}), (\ref{d3})
\item such that $|\we(A)-w(A)|<\delta $ for any $A\in V(H)$, where
the exact value of $\delta$ will be given later.
\end{itemize}
Now let $N$ be a natural number such that
\begin{itemize}
\item $\frac{N\we(A)}{l_A}\in\N$ if $A\in V(H)$.
\item $N\we(A,L_A,B)\in\N$ if $(A,L_A,B)\in E(H)$.
\end{itemize}

\vskip0.2in
\noindent
{\bf Step 1.} We construct an edge-less graph $Q$ such that:
\begin{itemize}
\item $V(Q)=\cup_{A\in V(H)} Q(A)$\quad (disjoint union)
\item $|Q(A)|=N \we(A)\,$
\item each $Q(A)$ is partitioned into
$\cup_{(A,L_A,B)\in E(H)} Q(A,L_A,B)$ such that $|Q(A,L_A,B)|=
N\we(A,L_A,B)$.
\end{itemize}
Since $\we$ satisfy our equations such $Q$ can be constructed.
\vskip0.2in
\noindent
{\bf Step 2.} We add edges to $Q$ in order to obtain the
graph $R$. For each pair $A,B$ that are connected in the graph $H$
form a bijection $Z_{A,B}:Q(A,L_A,B)\to Q(B,L_B,A)$.
If there is a loop in $H$ consider a bijection $Z_{A,A}$.
Then draw an edge between $x\in Q(A,L_A,B)$ and $y\in Q(B,L_B,A)$
if $Z_{A,B}(x)=y$.

\vskip0.2in
\noindent
{\bf Step 3.} Now we construct our graph $G$.
If $M\in U^{r+1,d}_n$ is a rooted labeled tree
such that $\mu_n(M)\neq 0$ let $Q(M)=\cup_{A\in R(M)} Q(A)$. We partition
$Q(M)$ into $\cup^{s_M}_{i=1} Q_i(M)$ such a way that each
$Q_i(M)$ contains exactly $l_A$ elements from the set $Q(A)$.
By the definition of $N$, we can make such partition.

\noindent
The elements of $V(G)$ will be the sets $\{Q_i(M)\}_{M\in 
 U^{r+1,d}_n\,,1\leq i \leq s_M}$. We draw one edge between
$Q_i(M)$ and $Q_j(M')$ if there exists $x\in Q_i(M), y\in Q_j(M')$
such that $x$ and $y$ are connected in $R$. We label the
vertex $Q_i(M)$ by the label of the root of $M$.
Let $Q_i(M)$ be a vertex of $G$ such that $M$ is a label-separated
tree. Note that if $M$ is not a rooted tree then $\mu_n(M)=0$.
 It is easy to see that the $r+1$-ball around $Q_i(M)$ in the graph $G$
is isomorphic to $M$ as rooted labeled balls. 
Also if $M$ is not label-separated then the $r+1$-ball around $Q_i(M)$
can not be a label-separated tree.
Therefore \begin{align} \label{becs2}
\sum_{L\in\urd_n\,,\mbox{ $L$ is not a label-separated tree}}
 p_G(L)
=\\=\sum_{L\in\urd_n\,,\mbox{ $L$ is not a label-separated tree}}
\sum_{A\in R(L)}\we(L)\leq \frac{\epsilon} {10} +\delta d |\urd_n|\,.
\end{align}
Also, if $M$ is a label-separated tree then
\begin{equation} \label{becs3}
|p_G(M)-\mu_n(T(\Grd,M))|\leq |R(M)|\delta\leq d\delta\,.
\end{equation}
Thus by (\ref{becs1}),(\ref{becs2}),(\ref{becs3})
if $\delta$ is choosen small enough then for any $\alpha\in U^{r+1,d}$
$$|p_G(\alpha)-\mu(T(\Grd,\alpha))|<\epsilon\,.$$
Thus our Theorem follows. \qed

\end{document}